\input amstex
\input amsppt.sty
\NoBlackBoxes
\nologo
\topmatter
\title Erratum to the paper ``A note on the factorization theorem of toric
birational maps after Morelli and its toroidal extension"
\endtitle
\author Kenji Matsuki \endauthor
\rightheadtext{Erratum}
\leftheadtext{Matsuki}

\endtopmatter

\document

This erratum describes:

\vskip.1in

1. the failure of the algorithm in [AMR] and [Morelli1] for the strong factorization
pointed out by Kalle Karu,

2. the statement of a refined weak factorization theorem for toroidal
birational morphisms in [AMR], in the form utilized in [AKMR] for the proof of the
weak factorization theorem for general birationla maps, avoiding the use of the
above mentioned algorithm for the strong factorization, and

3. a list of corrections for a few other mistakes in [AMR], mostly pointed out by
Laurent Bonavero.

We would like to emphasize that though [AMR] is a joint work with D. Abramovich
and S. Rashid, the author of this erratum, Kenji Matsuki, is solely responsible
for all the errors above.

\vskip.2in

\noindent 1. K. Karu pointed out that the procedure in Proposition 7.8 in [AMR]
does not preserve the condition $(*)$, contrary to its assertion, and thus the
proof does not work.  Moreover, the entire algorithm for the strong factorization
in Section 7 of [AMR] based upon Proposition 7.8, which attempted to correct the
logic of the line of argument of the original one described in [Morelli1] but is
identical to it as an algorithm, does not work, as shown by the following example of
a cobordism representing a toric birational map in dimension 3:

Consider a 4-dimensional simplicial cobordism $\Sigma$ describing three smooth
star subdivisions of the 3-dimensional cone $\langle v_1,v_2,v_3 \rangle$ at
$\langle v_1+v_2 \rangle$, $\langle v_2+v_3 \rangle$, and
$\langle v_1+v_2+v_3 \rangle$, in this order. The four maximal cones
in this cobordism are all pointing up, but the new cobordism $\Sigma'$, obtained by
the procedure of Proposition 7.8 in [AMR], contains a cone that is not pointing up.
Indeed, each of the four given maximal cones has exactly one positive extremal ray
and one extremal ray in the link of its circuit. The algorithm subdivides the
cobordism at the barycenters of the 2-dimensional cones generated by the positive
and the link extremal rays. After subdividing the two topmost cones, one of the new
cones will be a pointing up cone with one positive extremal ray $\langle
\rho\rangle$ and three negative extremal rays $\langle \rho_{12}\rangle$, $\langle
\rho_{23}
\rangle$ and $\langle \rho_3 \rangle$ with $n(\pi(\rho_{12})) = v_1 + v_2,
n(\pi(\rho_{23})) = v_2 + v_3$ and $n(\pi(\rho_3)) = v_3$. The next two subdivisions
are at the midrays (barycenters)
$\zeta_1 = \roman{Mid}(\langle
\rho_{12},\rho_{23}\rangle,l_{r(\langle
\rho_{12},\rho_{23}\rangle)})$ and $\zeta_2 = \roman{Mid}(\langle
\rho_{12},\rho_3\rangle,l_{r(\langle
\rho_{12},\rho_3\rangle)})$. The resulting subdivision $\Sigma'$ contains a cone
$\langle \rho, \rho_{12}, \zeta_1, \zeta_2 \rangle$, which is not pointing up as
it contains two positive extremal rays $\rho,
\rho_{12}$ and two negative extremal rays $\zeta_1,\zeta_2$.

\vskip.1in

\noindent 2. Weak Factorization Theorem for Toroidal Birational
Morphisms: \it Every proper and toroidal birational morphism $f:(U_X \subset X)
\rightarrow (U_Y \subset Y)$ between two nonsingular toroidal embeddings without
self-intersections can be factored into a sequence of blowups and blowdowns with
centers being smooth irreducible closed strata
$$\align
(U_X \subset X) &= (U_{V_0} \subset V_0) \dashrightarrow (U_{V_1} \subset V_1)
\dashrightarrow
\cdot\cdot\cdot \\
&\dashrightarrow (U_{V_{i-1}} \subset V_{i-1}) \dashrightarrow (U_{V_i} \subset
V_i)
\dashrightarrow \\
&\cdot\cdot\cdot
\dashrightarrow (U_{V_{l-1}} \subset V_{l-1})
\dashrightarrow (U_{V_l} \subset V_l) = (U_Y \subset Y) \\
\endalign$$
where all the toroidal embeddings are toroidal over $(U_Y \subset Y)$ and where
there is an index
$i_0$ such that for all
$i
\leq i_0$ the map
$V_i
\dashrightarrow X$ is a projective morphism and for all $i \geq i_0$
the morphism $V_i \rightarrow Y$ is a projective morphism.  In
particular, if both $X$ and $Y$ are projective then all
the $V_i$ are projective. 

\rm A proof of this theorem is given in the revised version (dated May 2000) of
[AKMW] via the process of
$\pi$-desingularization described in [AMR] (cf. [Morelli1,2]).  We remark that
this statement without the projectivity claim follows from our original argument
on Page 535 of [AMR].  Indeed, one only needs to take a simplicial,
$\pi$-nonsingular and collapsible cobordism
$\Sigma$ (in $(N_{\Delta})_{\Bbb Q}^+$ between $\Delta$ and $\Delta'$) which maps,
by definition, to
$\Delta$.  This can be obtained in the construction by requiring that the simplicial
and collapsible cobordism (in $N_{\Bbb Q}^+$ between
$\widehat{(\Delta_B^T)}$ and $\widehat{({\tilde \Delta}_B^T})$) should map to
$\widehat{(\Delta_B^T)}$, a requirement which is easily satisfied by first taking
the intersection of $\pi^{-1}(\widehat{(\Delta_B^T)})$ with any cobordism and
then further star subdividing if necessary. 

\vskip.1in

\noindent 3. List of corrections for some other mistakes

\vskip.1in

$\circ$ Example 4.1: We claimed that the presented cobordism $\Sigma$ in [AMR] was
$\pi$-nonsingular but it is actually $\pi$-singular.  An example of a
$\pi$-nonsingular and non-collapsible cobordsim, communicated to us by
Dan Abramovich, can be given as follows: We take vectors $v_1 = (1,0), v_2 = (0,1),
v_3 = (-1,-1) \in N_{\Bbb Q} = {\Bbb Z}^2
\otimes {\Bbb Q}$ and let $\Delta$ be the fan spanned by these vectors and hence
corresponding to ${\Bbb P}^2$ as a toric variety.  If we take the simplicial fan
$\Sigma$ in
$N_{\Bbb Q}^+ = N_{\Bbb Q} \oplus {\Bbb Q}$ whose maximal cones consist of
$$\align
\sigma_{11'2} &= \langle (v_1,0),(v_1,1),(v_2,0) \rangle \\
\sigma_{22'3} &= \langle (v_2,0),(v_2,1),(v_3,0) \rangle \\
\sigma_{33'1} &= \langle (v_3,0),(v_3,1),(v_1,0) \rangle \\
\sigma_{1'22'} &= \langle (v_1,1),(v_2,0),(v_2,1) \rangle \\
\sigma_{2'33'} &= \langle (v_2,1),(v_3,0),(v_3,1) \rangle \\
\sigma_{3'11'} &= \langle (v_3,1),(v_1,0),(v_1,1) \rangle, \\
\endalign$$
then $\Sigma$ is a simplicial $\pi$-nonsingular cobordism between $\Delta$ and
$\Delta$ itself, but is not collapsible as we have a directed cycle of the three
circuits
$$\sigma_{11'} = \langle (v_1,0),(v_1,1) \rangle \rightarrow \sigma_{22'} =
\langle (v_2,0),(v_2,1) \rangle \rightarrow \sigma_{33'} = \langle (v_3,0),(v_3,1)
\rangle \rightarrow \sigma_{11'}.$$

\vskip.1in

$\circ$ Definition 5.3:  We have to allow the equalities in the last line of the
definition \linebreak
`` ... in the set $\{\rho_i;r_i \leq 0\}$ or in the set
$\{\rho_i;r_i
\geq 0\}$, exclusively."

\vskip.1in

$\circ$ The definitions of $N_{\Delta}$ and
$(N_{\Delta})_{\Bbb Q}$ on Page 533 in [AMR] should be corrected as follows:

``Given a conical complex $\Delta$, we consider the space $(N^S \cap \sigma^S)
\oplus {\Bbb Z}$, for each $N^S = N_{\sigma^S}$ associated to the cone $\sigma^S
\in \Delta$, which can be glued together naturally via the glueing of $N^S \cap
\sigma^S$ to form the integral structure $N_{\Delta}$.  By considering the space
$\{(N^S \cap \sigma^S)
\oplus {\Bbb Z}\} \otimes {\Bbb Q}$ and glueing them together, we obtain the space
$$(N_{\Delta})_{\Bbb Q}^+ = (N_{\Delta} \oplus {\Bbb Z}) \otimes {\Bbb Q} =
(N_{\Delta})_{\Bbb Q} \oplus {\Bbb Q}."$$

We remark that a cobordism $\Sigma$ in $(N_{\Delta})_{\Bbb Q}^+$ between
$\Delta'$ and $\Delta$ for a refinement map $f_{\Delta}:\Delta' \rightarrow
\Delta$ of the conical complexes maps to $\Delta$ by definition.

$$\bold{References}$$

[AMR] D. Abramovich, K. Matsuki and S. Rashid, \it A note on the factorization
theorem of toric birational maps after Morelli and its toroidal extension, \rm
Tohoku Math. J. $\bold{51}$ (1999), 489-537

\vskip.1in

[AKMW] D. Abramovich, K. Karu, K. Matsuki and J. W{\l}odarczyk, \it Torification
and factorization of birational maps, \rm math.AG/9904135

\vskip.1in

[Morelli1] R. Morelli, \it The birational geometry of toric varieties, \rm J.
Alg. Geom. $\bold{5}$ (1996), 751-782

\vskip.1in

[Morelli2] R. Morelli, \it Correction to ``The birational geometry of toric
varieties", \rm homepage at the Univ. of Utah (1997), 767-770

\enddocument

\enddocument